\documentstyle[12pt]{article}   
\newcommand{\be}{\begin{equation}}
\newcommand{\ee}{\end{equation}}
\newcommand{\bea}{\begin{eqnarray}}
\newcommand{\eea}{\end{eqnarray}}
\newcommand{\bean}{\begin{eqnarray*}}
\newcommand{\eean}{\end{eqnarray*}}
\newcommand{\N}{I\!\!N}
\newcommand{\PP}{I\!\!P}
\newcommand{\R}{I\!\!R}
\newcommand{\Z}{Z\!\!\!Z}
\newcommand{\A}{A\!\!\!\!\!A}
\newcommand{\G}{I\!\!L}

\newcommand{\suml}{\sum\limits}

\newcommand{\C}{I\!\!\!\!C}
\newcommand{\Q}{I\!\!\!\!Q}

\newcommand{\Spec}{\mbox{Spec$\;$}}

\newcommand{\veee}{\scriptscriptstyle\vee}
\newcommand{\conv}{\mbox{conv$\;$}}

\newcommand{\linh}{\mbox{span}}

\newcounter{Abschnitt}[section]
\newcommand{\neu}[1]{\protect\refstepcounter{Abschnitt}\protect\label{t#1}\vspace{1ex}
                     {\bf (\arabic{section}.\protect\arabic{Abschnitt})}
                     $\qquad$}
\newcommand{\zitat}[2]{(\protect\ref{s#1}.\protect\ref{t#1#2})}
\newcommand{\surj}{\longrightarrow\hspace{-1.5em}\longrightarrow}

\newcounter{secnum}

\newcommand{\sect}[1]
 {\protect\section{#1}
  \protect\setcounter{secnum}{\value{section}}
  \protect\setcounter{equation}{0}
 \protect\renewcommand{\theequation}{\mbox{\arabic{secnum}.\arabic{equation}}}}

\begin{document}
\title{Toric $\,\Q$-Gorenstein Singularities}

\author{Klaus Altmann\\
	 \small Dept. of Mathematics, M.I.T., Cambridge, MA 02139, U.S.A.
	 \vspace{-0.7ex}\\ \small E-mail: altmann@math.mit.edu}
\date{}
\maketitle

\begin{abstract}
For an affine, toric $\,\Q$-Gorenstein variety $Y$ (given by
a lattice polytope $Q$) the vector space $T^1$ of infinitesimal deformations
is related to the complexified vector spaces
of rational Minkowski summands of faces of $Q$.\\
Moreover,
assuming $Y$ to be an isolated, at least 3-dimensional singularity, $Y$ will
be rigid unless it is even Gorenstein and $\mbox{dim}\,Y=3\;
(\mbox{dim}\,Q=2)$.\\
For this particular case,
so-called toric deformations of $Y$ correspond to Min\-kow\-ski decompositions
of $Q$ into a sum of lattice polygons. Their Kodaira-Spencer-map can be
interpreted in a very natural way.\\
We regard the projective variety $\PP (Y)$ defined by the lattice
polygon $Q$. Data concerning the deformation theory of $Y$ can be interpreted
as data concerning the Picard group of $\PP (Y)$.\\
Finally, we provide some examples (the cones over the toric Del Pezzo
surfaces). There is one such variety yielding $\Spec \C[\varepsilon]/_
{\displaystyle \varepsilon^2}$ as the base space of the semi-universal
deformation.
\end{abstract}

\tableofcontents
\par
\vspace{2ex}



%
%
\sect{Introduction}\label{s1}


\neu{11}
In \cite{T1} and \cite{Mink} we investigated the deformation theory of
affine toric varieties $Y=\Spec \C [\stackrel{\veee}{\sigma}\cap M]$
(cf. \zitat{2}{1} for an explanation of the notations):\\
\par
It is always the first step to look at the vector space $T^1_Y$ of
infinitesimal deformations - it equals (if $Y$ admits isolated singularities
only) the tangent space on the base space $S$ of the semi-universal
deformation of $Y$. For toric $Y$ the space $T^1_Y$ is $M$-graded, and the
homogeneous pieces were computed in \cite{T1} (cf. \zitat{2}{2} of the present
paper).\\
\par
In \cite{Mink} we were interested in describing toric deformations of $Y$.
They are defined as those deformations (i.e. flat maps $f:X\rightarrow S$
endowed with an isomorphism $f^{-1}(0\in S) \cong Y$) such that the total
space $X$ together with the embedding of the special fiber $Y \hookrightarrow
X$
are toric.\\
\par
Toric deformations are ``really existing'' deformations -
in the sense that they admit reduced (even smooth) base spaces. Moreover,
we conjecture that (in case of isolated singularities $Y$)
the semi-universal deformation of $Y$ is toric over each irreducible component
of the reduced base space. (This is true for $\mbox{dim}\, Y = 2$ and also for
examples of higher dimension.)\\
\par
Toric deformations always arise as relative deformations
of $Y$ inside a greater affine toric variety $X$ containing $Y$ as a relatively
complete intersection. More strictly speaking, $Y\subseteq X$ is defined by a
so-called toric regular sequence $x^{r_0^1}-x^{r_1^1},\dots,
x^{r_0^m}-x^{r_1^m} \in \Gamma(X, {\cal O}_X)$.\\
\par
On the other hand, each toric regular sequence can be regarded as a flat map
$X\rightarrow \C^m$
by itself. This $m$-parameter deformation of $Y$ is called the standard
toric deformation induced by the given sequence.\\
\par
It is possible to compute the Kodaira-Spencer-map $\varrho: \C^m \rightarrow
T^1_Y$ corresponding to these standard toric deformations. Then, the first
observation is that $\varrho$ maps
the $i$-th canonical basis vector $e^i\in \C^m$ into the homogeneous piece
$T^1_Y(-\bar{r}^i)$, and $\bar{r}^i \in M$ is defined as the exponent of the
common image of $x^{r^i_0}$ and $x^{r^i_1}$ via the surjection
$\Gamma(X, {\cal O}_X)\surj\Gamma(Y, {\cal O}_Y)$.\\
\par


{\bf Definition:}
A toric regular sequence is called strongly homogeneous if
only $m+1$ different elements occur in the set $\{r^1_0,\,r^1_1,\dots,
\,r^m_0,\,r^m_1\}$. Then, the corresponding images $\bar{r}^1,\dots,\bar{r}^m$
coincide, and this element will be denoted by $\bar{r}\in
\stackrel{\veee}{\sigma}\cap M$.
It equals the
negative degree of the Kodaira-Spencer-map.\\
(Toric regular sequences of length one are always strongly homogeneous.)\\
\par

The main result of \cite{Mink} is a complete combinatorial description of
those standard toric deformations that are induced by strongly homogeneous
toric regular sequences. They arise from certain Minkowski decompositions
of affine slices (induced by $\bar{r}\in M$) of the cone $\sigma$ (cf.
\zitat{3}{1} of the present paper).\\
\par


\neu{12}
The aim of this paper is to apply the previous results to the special case
of toric $\,\Q$-Gorenstein singularities. On the one hand, this notion is the
next one if we are looking for a wider class than that of complete
intersections (which yields no interesting deformation theory). On the other
hand, the property ``$\,\Q$-Gorenstein'' admits a very clear description
in the language of toric varieties and convex cones:\\
\par
In general, the dualizing sheaf $\omega$ on a Cohen-Macaulay variety
is defined as
\begin{itemize}
\item[(i)]
$\omega_P := \Omega_P^{\,\mbox{\scriptsize dim}\,P}$ (sheaf of the highest
differential forms) if $P$ is smooth, and
\item[(ii)]
$\pi_{\ast}\omega_Y:= \mbox{\underline{Hom}}_{{\cal O}_P}(\pi_{\ast}{\cal
O}_Y,\,
\omega_P)$ for flat and finite maps $\pi: Y \rightarrow P$.
\vspace{1ex}
\end{itemize}
\par
%
%


{\bf Definition:} A variety $Y$ is called ($\,\Q$-) Gorenstein if (the
reflexive hull of some tensor power of) $\omega_Y$ is an invertible sheaf
on $Y$.\\
\par

Since toric varieties are normal, the dualizing sheaf can be obtained as the
push forward of the canonical sheaf on its smooth part. Hence, in our
special situation, $\omega_Y$ equals the $T\mbox{(orus)}$-invariant complete
fractional ideal that is given by the order function mapping each fundamental
generator onto $1\in \Z$ (cf. Theorem I/9 in \cite{Ke}).\\
\par
In particular, we obtain the following\\
\par


{\bf Fact:} Let $Y=\Spec \C [\stackrel{\veee}{\sigma}\cap M]$ be an affine
toric variety given by a cone $\sigma= \langle a^1,\dots, a^N \rangle$. (The
fundamental generators $a^i$ are assumed to be primitive elements of the
lattice that is dual to $M$.)\\
Then, $Y$ is $\,\Q$-Gorenstein, if and only if
there is a primitive element $R^{\ast}\in M$ and a natural number
$g\in \N$ such that
\[
\langle a^i,\, R^{\ast}\rangle = g \quad \mbox{for each } i=1,\dots,N.
\]
$Y$ is Gorenstein if and only if $g=1$, in addition.\\
\par


\neu{13}
Affine toric varieties of dimension two are always $\,\Q$-Gorenstein. The
deformation theory of these varieties (the two-dimensional cyclic quotient
singularities) is well known. For instance, $T^1_Y$ and the
number and dimension of the components of the reduced base space of the
semi-universal deformation have been computed (cf. \cite{Riem}, \cite{Arndt},
\cite{Ch}, \cite{St}).\\
Therefore, our investigation concerns the case of $Y$ being smooth in
codimension 2.\\
\par
Assume that $Y=\Spec \C [\stackrel{\veee}{\sigma}\cap M]$ is a
$\,\Q$-Gorenstein
variety (i.e. $\langle a^i,\, R^{\ast}\rangle = g$ for each $i=1,\dots,N$),
which is smooth in codimension 2. Denote by $Q$ the lattice polyhedron
$Q:= \mbox{Conv}(a^1,\dots,a^N)$.
Then, we obtain the following results:
\begin{itemize}
\item[(1)]
The graded pieces of $T^1_Y$ equal the vector spaces of Minkowski summands
of certain faces of $Q$ (cf. Theorem \zitat{2}{7}).
\item[(2)]
For the special case of an isolated singularity $Y$ this implies:
\begin{itemize}
\item
If $g\geq 2$, then $Y$ will be rigid, i.e. $T^1_Y=0$.
\item
If $Y$ is Gorenstein (i.e. $g=1$) and at least 4-dimensional, then $Y$ will
be rigid.
\item
Let $Y$ be a 3-dimensional Gorenstein singularity given by a plane, convex
$N$-gon $Q$. Then, $T^1_Y$ is concentrated in the single degree $-R^\ast$,
and it equals to the $(N-3)$-dimensional vector space of Minkowski summands of
$Q$.
\end{itemize}
(Cf. \zitat{2}{8} and \zitat{2}{9}.)
\item[(3)]
Keep assuming that $Y$ is a 3-dimensional, isolated, toric, Gorenstein
singularity. Then,
toric $m$-parameter deformations of $Y$
correspond to Minkowski decompositions of $Q$ into a sum of $m+1$ lattice
polygons. Using the previous description of $T^1_Y$, the Kodaira-Spencer-map
is the natural one.
\vspace{2ex}
\end{itemize}
\par


\neu{14}
Let us assume, for a moment, that the semi-universal deformation is, indeed,
toric over each component of the reduced base space $S_{\mbox{\scriptsize
red}}$.
Then, the absence of proper lattice summands of $Q$
would mean that $S$ is only 0-dimensional.\\
\par
On the other hand, this is possible even for $Y$ admitting a non-trivial
$T^1_Y$ (which equals the tangent space of $S$).
In (\ref{s4}.\ref{t44}.3) we will give a (three-dimensional) example of this
phenomenon:
The base space $S$ is a point with non-reduced structure, i.e. each deformation
is obstructed.\\
\par


\neu{15}
A toric Gorenstein variety $Y=\Spec \C [\stackrel{\veee}{\sigma}\cap M]$
is given by a lattice polytope $Q$ ($\sigma = \mbox{Cone}(Q)$).
On the other hand, each lattice polytope $Q$ induces a projective toric
variety (defined by the inner normal fan) endowed with an ample line
bundle ${\cal O}(1)$. We call this the polar variety of $Y$ - denoted by
${\PP}(Y)$.\\
Then, the cone $P(Y)$ over the embedded projective variety
${\PP}(Y)$ equals the affine toric variety which is given by the cone dual to
that of $Y$.\\
\par
In \S \ref{s4}, we  interprete data of the deformation theory of
3-dimensional $Y$ as
data concerning divisors on the varieties ${\PP}(Y)$ or $P(Y)\setminus
\{0\}$. Assuming
$Y$ having only an isolated singularity
we obtain the following relations:
\begin{itemize}
\item[(4)]
$T^1_Y \,=\; \mbox{Pic}(P(Y)\setminus \{0\}) \otimes_{\Z} \C \,=\;
\mbox{Pic}\,{\PP}(Y) \left/_{\displaystyle \!\!{\cal O}(1)}\right.
\otimes_{\Z} \C$.
\item[(5)]
Toric $m$-parameter deformations of $Y$ correspond to splittings of
${\cal O}_{{\PP}(Y)}(1)$ into a tensor product
of $m+1$ invertible sheaves that are nef.
\vspace{2ex}
\end{itemize}


\neu{16}
{\bf Acknowledgement:}
I am grateful to Duco van Straten for computing several semi-universal
deformations on the computer (using Macaulay).\\
Moreover, I want to thank Bernd Sturmfels for his special lecture
concerning the cone of Minkowski summands and the Gale transform of a given
polytope.\\
\par

%
%
\sect{The $T^1$ of toric $\,\Q$-Gorenstein singularities}\label{s2}


\neu{21}
Let us start with introducing some basic notations and recalling the
$T^1$-formulas of \cite{T1}:\\
\par
Let $M, N$ be free $\Z$-modules of finite rank - endowed with a perfect
pairing $\langle .,.\rangle : N \times M \rightarrow \Z$. Denote by
$M_{\R}$ and $N_{\R}$ the corresponding vector spaces (dual to each other)
obtained via base change with $\R$.\\
\par
Let $\sigma= \langle a^1,\dots, a^N \rangle \subseteq N_{\R}$ be a
top-dimensional, rational, polyhedral cone with apex.
The fundamental generators $a^i \in N$ are assumed
to be primitive elements of the lattice $N$.\\
\par
The dual cone $\stackrel{\veee}{\sigma}\subseteq M_{\R}$ of $\sigma$ is
defined as $\stackrel{\veee}{\sigma}:= \{r\in M_{\R} |\; \langle a,\, r\rangle
\geq 0 \;\mbox{ for each } a\in \sigma\}$. Denote by $E \subseteq
\stackrel{\veee}{\sigma}\cap M$ the minimal (finite) set that
generates the semigroup $\stackrel{\veee}{\sigma}\cap M$.
In particular, $Y:= \mbox{Spec}\,\C [\stackrel{\veee}{\sigma}\cap M] \subseteq
\C^E$.\\
\par


\neu{22} {\bf Theorem:} (cf. (2.3) and (4.4) of \cite{T1})
\begin{itemize}
\item[(1)]
The vector space $T^1_Y$ of infinitesimal deformations of $Y$ is $M$-graded.
For a fixed element $R\in M$ the homogeneous piece of degree $-R$ can be
computed as
\[
T_Y^1(-R)=\left(^{\displaystyle
L(E')}\!\left/{\displaystyle\suml^N_{i=1}L(E_i)}
^{\makebox[0mm]{}}\right.\right)^{\ast}\otimes_{\R} \C
\]
$(E_i := \{s\in E\;|\;\;
(0\le)\langle a^i,s\rangle<\langle a^i,R\rangle\}
\quad (i=1,\ldots,N)\,;\quad
E' :=\; \bigcup\limits^N_{i=1} E_i\,;$\\
$L(\mbox{set}) := \mbox{$\R$-vector space of all linear dependences between
its elements}).$
\item[(2)]
If $Y$ is smooth in codimension 2 (i.e. if all 2-dimensional faces
$\langle a^i, a^j \rangle < \sigma$ are spanned by a part of a $\Z$-basis of
the lattice $N$), then
with
\[
V_i := {\linh}_{\R}(E_i)=\left\{
\begin{array}{l@{\quad\mbox{for}\;\:}l}
0 &\langle a^i,R\rangle\le 0\\
\left[ a^i=0\right] \subseteq M_{\R} & \langle a^i,R\rangle =1\\
M_{\R} & \langle a^i,R\rangle\ge2
\end{array}
\right.
\]
we obtain the second formula
\[
T_Y^1(-R)=\mbox{Ker}\left[^{\displaystyle V_1\oplus\ldots\oplus V_N}\left/
_{\!\!\!\displaystyle\sum_{\langle a^i,a^j\rangle <\sigma}V_i\cap V_j}\right.
\surj
V_1+\ldots+V_N\right]^{\ast}.
\]
\end{itemize}
\par


\neu{23} {\bf Lemma:}
Let $Y$ be a $\,\Q$-Gorenstein variety, which is smooth in codimension 2.
If $R\in M$
is a degree such that $\langle a^i,\,R \rangle \geq 2$ for some $i \in
\{1,\dots,N\}$, then $T^1_Y(-R)=0$.\\
\par


{\bf Proof:}
Let $R^{\ast}\in M$ be as in \zitat{1}{2}, i.e. $\langle a^i,\,R^{\ast}
\rangle = g$ for $i=1,\dots,N$. Then,
\[
H:= \{a\in N_{\R}|\; \langle a,\, g\,R-R^{\ast}\rangle = 0\}
\]
is a hyperplane in $N_{\R}$ that subdivides the set of fundamental generators
of $\sigma$. $H^-, H$, and $H^+$ contain the elements $a^i$ meeting the
properties $\langle a^i,\,R\rangle \leq 0$, $\langle a^i,\,R\rangle =1$, and
$\langle a^i,\,R\rangle \geq 2$, respectively. Let us assume that the latter
class of generators is not empty.\\
\par
To use the second $T^1$-formula of the previous theorem, we fix a map
\[
\varphi: \{i\,|\; \langle a^i,R\rangle =1\} \rightarrow \{1,\dots,N\}
\]
such
that for each $a^i\in H$ the element $a^{\varphi(i)}$ is contained in $H^+$
and adjacent to $a^i$ (i.e. $\{a^i,\, a^{\varphi(i)}\}\subseteq\sigma$ spans a
2-dimensional face
of $\sigma$).\\
\par
Now, assume that we are given an element $v=(v_1,\dots,v_N)\in V_1\oplus
\dots\oplus V_N$ such that $v_1+\dots +v_N=0$. Adding the terms
$[-v_i\cdot e^i+v_i\cdot e^{\varphi(i)}]$ (for $\langle a^i,R\rangle =1$)
does not change the equivalence class of $v$ in
$^{\displaystyle V_1\oplus\ldots\oplus V_N}\!\!\left/
_{\!\!\!\!\!\displaystyle\sum_{\langle a^i,a^j\rangle <\sigma}
V_i\cap V_j}\right.$.
However,
non-trivial components survive for $\langle a^i,R\rangle \geq 2$ (corresponding
to $V_i=M_{\R}$) only.\\
\par
The set of these special generators $a^i$ is connected  by 2-dimensional
faces of $\sigma$. Moreover, by slightly disturbing $R$ inside $M_{\R}$, we
can find a unique $a^{\ast}$ among these edges on which $R$ is maximal. Then,
each $a^i\in H^+$ is connected with $a^{\ast}$ by an $R$-monotone path
(consisting of 2-faces of $\sigma$) inside $H^+$.\\
Now, we can use the previous method of cleaning the components of $v$ once
more - the steps from $a^i$ to $a^{\varphi(i)}$ are replaced by the steps on
the path from $a^i$ to $a^{\ast}$. It remains an $N$-tuple $v$ which is
non-trivial at most at the $a^{\ast}$-place. On the other hand, the components
of $v$ sum up to $0$, but this yields $v=0$.
\hfill$\Box$\\
\par


\neu{24}
If $\langle a^i,\, R\rangle \leq 1$ for every $i \in \{1,\dots,N\}$, then
equality holds on some face $\tau < \sigma$. Now, $\tau$ is a
top-dimensional cone in the linear subspace $\tau - \tau \subseteq N_{\R}$,
and it defines a variety $Y_{\tau} = \mbox{Spec} \;\C[
^{\displaystyle\tau^{\veee}\cap M}\!/\!_{
\displaystyle\tau^{\bot}\cap M} ]$, which is even Gorenstein.
The corresponding element $R^{\ast}_{\tau} \in
^{\displaystyle M}\!/\!_{\displaystyle \tau^{\bot}\cap M}$
can be obtained as the image of $R$ as well as of $\frac{1}{g}R^{\ast}$
using the canonical projection $M_{\R} \surj
^{\displaystyle M_{\R}}\!/\!_{\displaystyle \tau^{\bot}}$.
\\
\par


{\bf Lemma:}
In general (even the $\,\Q$-Gorenstein assumption can be dropped)
let $\tau < \sigma$ be a face such that $\langle a^i, R\rangle
\geq 1$ for $a^i\in\tau$ and $\langle a^i, R\rangle \leq 0$ otherwise. Then,
$T^1_Y(-R) = T^1_{Y_{\tau}}(-[\mbox{image of }R])$.\\
In the special situation discussed previously this means
$T^1_Y(-R) = T^1_{Y_{\tau}}(-R^{\ast}_{\tau})$.\\
\par
{\bf Proof:}
The formula of Theorem \zitat{2}{2}(1) remains true if $E$ is replaced by
an arbitrary (not necessarily minimal) generating subset of $\sigma^{\veee}
\cap M$ - even a multiset could be allowed. Hence, for computing
$T^1_{Y_{\tau}}(-R^{\ast}_{\tau})$ we can use the image $\bar{E}$ of $E$
under the projection
\[
\sigma^{\veee}\cap M \surj ^{
\displaystyle (\sigma^{\veee}+\tau^{\bot})\cap M}\!\! \left/
\!\! _{\displaystyle \tau^{\bot}\cap M} \right. .
\]
For $a^i\in \tau$, the corresponding sets $\bar{E}_i \subseteq \bar{E}$
coincide with the images of the subsets $E_i$. For $a^i \notin \tau$, the
notion $ \bar{E}_i$ does not make sense, and the $E_i$ are empty, anyway. It
remains to show that the canonical map
\[
\left. ^{\displaystyle L(E')} \!\! \right/ \!\!
_{\displaystyle \sum_{a^i\in \tau} L(E_i)} \longrightarrow
\left. ^{\displaystyle L(\bar{E}')} \!\! \right/ \!\!
_{\displaystyle \sum_{a^i\in \tau} L(\bar{E}_i)}
\]
is an isomorphism.\\
\par
The vector space $\tau^{\bot}$ is generated by
$\tau^\bot\cap (\bigcap_{a^i\in\tau} E_i)$.
Hence, by choosing a basis among these elements, we can embed $\tau^{\bot}$
into $\R^{\tau^\bot\cap (\cap_\tau E_i)}$ to obtain a section of
\bean
\R^{\tau^\bot\cap (\cap_\tau E_i)} \subseteq
\bigcap_{a^i\in\tau} L(\bar{E}_i) \subseteq
L(\bar{E}') &\surj& \tau^\bot \\
(\dots,\lambda_r,\dots)_{r\in E'} & \mapsto & \sum_{r\in E'} \lambda_r \cdot
r \in \tau^\bot \subseteq M_{\R}.
\eean
In particular, we obtain
\bean
L(\bar{E}_i) &=& L(E_i) \oplus \tau^\bot \;\; (a^i \in \tau)
\qquad\mbox{ and}\\
L(\bar{E}') &=& L(E') \oplus \tau^\bot.
\qquad\qquad\qquad\qquad\qquad\qquad\qquad\qquad\qquad\Box
\eean
\par


\neu{25}
Let $Q=\mbox{Conv}(a^1,\dots,a^N)$ be a $K$-polytope contained in some
$K$-vector space $\A$ ($K= \Q\mbox{ or } \R$). It can be described by
inequalities
\[
\langle \bullet,\, -c^v\rangle \leq \eta^v \quad (c^v\in \A^{\ast},\,
\eta^v \in K)
\]
corresponding to the facets of $Q$.
(The $c^v$ are the inner normal vectors of $Q$.)\\
\par
Let us denote by $\sigma \subset \A\times K$ the cone over $Q$ (embedded
as $Q\times\{1\}$). Then, the
pairs $[c^v,\eta^v]\in \A^{\ast}\times K$ are exactly the fundamental
generators of the dual cone $\stackrel{\veee}{\sigma}$.
For $i=1,\dots,N$ we define
\bean
F_i &:=& \{[c^v,\eta^v]\in \A^{\ast}\times K \,|\; \langle a^i,\,-c^v
\rangle = \eta^v\}\\
&=& \{\mbox{fundamental generators of the face } (a^i)^{\bot}\cap
       \stackrel{\veee}{\sigma} \,<\; \stackrel{\veee}{\sigma}\}.
\eean
Finally, denoting the set of all
fundamental
generators of $\stackrel{\veee}{\sigma}$
by $F':=\; \bigcup\limits^N_{i=1} F_i$,
we can construct the
$K$-vector space
\[
\tilde{T}^1(Q):=
\left(^{\displaystyle L(F')}\!\left/{\displaystyle\suml^
N_{i=1}L(F_i)}^{\makebox[0mm]{}}\right.\right)^{\ast}.
\]
Now, each Minkowski summand $Q'$ of a scalar multiple of $Q$ is given by
inequalities $\langle \bullet,\, -c^v\rangle \leq {\eta '\,}^v$. We can define
its class $\varrho(Q')\in \tilde{T}^1(Q)$ as
\[
\varrho(Q')\, (q\in L(F')):= \sum_v q_v {\eta '\,}^v \in K.
\]
This definition is correct and depends on the translation class of $Q'$
only. Moreover, scalar multiples of $Q$ yield the zero class.\\
\par
On the other hand, the translation classes of Minkowski summands of scalar
multiples of $Q$ form a convex polyhedral cone which contains ``$Q$'' as an
interior point. Dividing by the relation ``$Q=0$'' yields a $K$-vector space
which we will call the vector space of Minkowski summands of scalar multiples
of $Q$. It has one dimension less than the cone of Minkowski summands.\\
\par


{\bf Theorem:} (cf. \cite{Sm})\quad
The map $\varrho$ induces an isomorphism between the vector space of
Minkowski summands of scalar multiples of $Q$ and the vector space
$\tilde{T}^1(Q)$.\\
\par


{\bf Remark:} The constructions made in \zitat{2}{5} do not depend on the
linear, but on the affine structure of $\A$.\\
\par


\neu{26}
{\bf Lemma:}
Let $Y$ be an affine
toric Gorenstein variety induced
from a lattice
polytope $Q$. Then, the vector space $T^1_Y(-R^\ast)$
equals the complexified vector space $\tilde{T}^1(Q)\otimes \C$ of
(rational or real) Min\-kow\-ski summands of scalar
multiples of $Q$ (modulo translations and scalar multiples of $Q$ itself).\\
\par


{\bf Proof:}
We will use the first $T^1$-formula of Theorem \zitat{2}{2}. For the special
degree $-R^{\ast}$ the sets $E_i$ equal
\[
E_i= \{s\in E\,|\; \langle a^i,\, s\rangle =0\} = E\cap (a^i)^{\bot}.
\]
In particular, they contain the sets $F_i$ constructed above. We obtain a
natural linear map
\[
\theta: T^1_Y(-R^\ast) = \left(^{\displaystyle L(E')}\!\left/
{\displaystyle\suml^
N_{i=1}L(E_i)}^{\makebox[0mm]{}}\right.\right)^{\ast}\otimes_{\R} \C
\longrightarrow
\left(^{\displaystyle L(F')}\!\left/{\displaystyle\suml^
N_{i=1}L(F_i)}^{\makebox[0mm]{}}\right.\right)^{\ast}\otimes_{\R} \C
= \tilde{T}^1(Q)\otimes \C,
\]
and it remains to prove that $\theta$ is an isomorphism.\\
\par
Let $s\in E'\subseteq \partial\sigma^{\veee}$ be an element that
is not a fundamental generator of $\stackrel{\veee}{\sigma}$. Then, there
is a minimal face $\alpha < \sigma^{\veee}$ containing $s$ (as a
relatively interior point), and we can choose some fundamental generators
$s^1,\dots, s^k \in \alpha$ such that $s=\sum_{j=1}^k {\lambda}_j\,s^j\;
({\lambda}_j \in \Q_{\geq 0})$.\\
\par
Now, for each $E_i$ containing $s$ (equivalent to $\alpha < (a^i)^{\bot}
\cap \stackrel{\veee}{\sigma}$) we have a decomposition
\[
L(E_i) = L(E_i \setminus \{s\}) \;\oplus\; \R\cdot [\mbox{relation }
s=\sum_{j=1}^k {\lambda}_j\,s^j].
\]
In particular, the second summand can be reduced in the expression for
$T^1_Y$. Since the map $\theta$ consists of such steps only, we are done.
\hfill$\Box$\\
\par


\neu{27}
We collect the results of \zitat{2}{3} - \zitat{2}{6}: Let $Y$ be an affine
toric $\,\Q$-Gorenstein variety given by a lattice polytope $Q =
\mbox{Conv}(a^1,\dots,a^N)$ contained in an affine hyperplane $[\langle
\bullet , R^\ast\rangle = g] \subseteq N_{\R}$ of lattice-distance $g$ from
$0\in N_{\R}$. Moreover, assume that $Y$ is smooth in codimension 2.\\
Then, the graded pieces of $T^1_Y$ are related to the vector spaces of
Minkowski summands of faces of $Q$. Using the notations of \zitat{2}{5} (and
$\tilde{T}^1(\emptyset):= 0$) we obtain the following two equivalent
descriptions of $T^1_Y$:\\
\par
{\bf Theorem:}
\begin{itemize}
\item[(1)]
Let $R\in M$, then
\bean
T^1_Y(-R) &=&
\left\{ \begin{array}{ll}
\tilde{T}^1(Q\cap[\frac{1}{g}R^\ast-R]^\bot) &=
\tilde{T}^1(\mbox{Conv}\{a^i\,|\,\langle a^i,R\rangle =1\})\\
& \mbox{ for } \frac{1}{g}R^\ast \geq R \mbox{ on } \sigma \;(\mbox{i.e.}
\langle a^i,R\rangle\leq 1 \;\forall i)\\
0 & \mbox{ otherwise}.
\end{array} \right.
\eean
\item[(2)]
Let $\tau < \sigma$ be a face of $\sigma$, then
$T^1_Y\left( [-\frac{1}{g}R^\ast +
\mbox{int}(\sigma^{\veee}\cap\tau^\bot)]\cap M \right) = \tilde{T}^1
(Q\cap \tau)$. $T^1_Y$ vanishes in the remaining degrees.
\vspace{2ex}
\end{itemize}
\par


\neu{28}
With the same assumptions as in \zitat{2}{7} we obtain the following
applications of the previous theorem:
\pagebreak
\\
\par
{\bf Corollary:}
\begin{itemize}
\item[(1)]
If every 2-face of $Q$ is a triangle (for instance, if $Y$ is smooth in
codimension 3), then $Y$ is rigid, i.e. $T^1_Y=0$.
\item[(2)]
If $Y$ is Gorenstein ($g=1$) of dimension at least 4 ($\mbox{dim}\,Q\geq 3$),
then the existence of a 2-face of $Q$ that is not a triangle implies
$\mbox{dim}\, T^1_Y = \infty .$
\item[(3)]
Let $Y$ be not Gorenstein, i.e. $g\geq 2$. Then, $\mbox{dim}\, T^1_Y <
\infty$ implies $T^1_Y=0$.
\vspace{2ex}
\end{itemize}
\par

{\bf Proof:}
(1) If $Q$ is shaped that every 2-face is a triangle, then every (at least
2-dimensional) face of $Q$ will have this property, too. On the other hand,
Smilanski has shown that polytopes with only triangular 2-faces admit at
most trivial Minkowski decompositions (cf. \cite{Sm}, Corollary (5.2)), i.e.
$\tilde{T}^1=0$.\\
Moreover, faces of dimension smaller or equal than 1 of $Q$ cannot be
non-trivially decomposed, anyway.\\
\par
(2) Two-dimensional polygons with at least 4 vertices have a non-trivial
$\tilde{T}^1$. Hence, a non-triangular 2-face of $Q$ yields a proper face
$\tau < \sigma$ with $\tilde{T}^1(Q\cap \tau) \neq 0$.\\ On the other hand,
``proper'' means that $[-R^\ast + \mbox{int}(\sigma^{\veee}\cap\tau^\bot)]
\cap M$ contains infinitely many elements, and $T^1_Y$ is non-trivial in all
those degrees.\\
\par
(3) Assume that $T^1_Y \neq 0$, then there must be a face $\tau < \sigma$
and an element $-R\in [-\frac{1}{g}R^\ast +
\mbox{int}(\sigma^{\veee}\cap\tau^\bot)]\cap M$ such that $T^1_Y(-R) =
\tilde{T}(Q\cap\tau) \neq 0$.\\
For $\tau = \sigma$ we would obtain
$[-\frac{1}{g}R^\ast +
\mbox{int}(\sigma^{\veee}\cap\tau^\bot)]\cap M = \{-\frac{1}{g}R^\ast\}
\cap M = \emptyset$. Hence, $\tau<\sigma$ must be a proper face, and we can
argue as in (2).
\hfill$\Box$\\
\par


\neu{29}
Finally, we want to mention the case $\mbox{dim}\,Y = 3$. Then, $Y$ is an
isolated singularity, and it is given by a lattice $N$-gon $Q$.
\\
\par
{\em Case1: $Y$ is not Gorenstein (i.e. $g\geq 2$).}\\
Then, $Y$ is rigid. (This follows from (3) of the corollary in
\zitat{2}{8}.)\\
\par
{\em Case 2: $Y$ is Gorenstein (i.e. $g=1$).}\\
Then, $T^1_Y$ is concentrated in degree $-R^\ast$, and $T^1_Y=T^1_Y(-R^\ast)
= \tilde{T}^1(Q)$ has dimension $N-3$.\\
(The proper faces of $Q$ are Minkowski indecomposable. Hence, to produce a
non-trivial contribution to $T^1_Y$, the face $\tau$ in Theorem
\zitat{2}{7}(2) has to equal $\sigma$. That means, $T^1_Y$ is concentrated
in degree $-R^\ast$ only.\\
On the other hand, for computing the dimension of $T^1_Y(-R^\ast)
= \tilde{T}^1(Q)$ in our special case, use the second formula of Theorem
\zitat{2}{2}.)\\
\par

%
%
\sect{Really existing deformations of toric Gorenstein singularities}\label{s3}


\neu{31}
As in the previous chapter, we start with recalling the general result
concerning arbitrary affine toric varieties. We use the notations of
\zitat{2}{1}.\\
\par

Let $\bar{r}\in \stackrel{\veee}{\sigma}\cap M$ be a primitive element.
Then, each (strongly homogeneous)
toric regular sequence of degree $-\bar{r}$ and its corresponding standard
deformation of $Y$ arise in the following way:
\begin{itemize}
\item[(i)]
Define $(\A,\,\G)$ as the affine space (with lattice) induced by $\bar{r}\in
M$
\[
(\A,\,\G) := (N_{\R},\,N) \cap \{a\in N_{\R}|\; \langle a,\,\bar{r} \rangle
= 1 \}.
\]
By choosing an arbitrary base point $0\in \G$ the pair $(\A,\,\G)$ can be
regarded as a vector space with lattice. Moreover, we obtain an isomorphism
of lattices $\G\times\Z \stackrel{\sim}{\rightarrow} N$ via
$(a,\,g) \mapsto (a-0)+g\cdot 0$.
\item[(ii)]
$Q:= \sigma \cap \A$ is a (not necessarily compact) rational polyhedron
in $\A$.
Fix a Minkowski decomposition $Q=R_0 + \dots + R_m$ such that for each vertex
of $Q$ at least $m$ of its $m+1$ $R_i$-summands (which are uniquely determined
vertices of the polyhedra $R_i$) are contained in the lattice $\G$.
\item[(iii)]
Define $P\subseteq \A\times \R^{m+1}$ as the convex polyhedron
\[
P:= \conv \left( \bigcup_{i=0}^m R_i\times \{e^i\} \right)
\]
and $\tilde{\sigma}:= \overline{\R_{\geq 0} \cdot P}$ as the closure of its
cone in $\A\times\R^{m+1}$.\\
Moreover, if $r^i$ denotes the projection of $\G\times\Z^{m+1}$ onto the $i$-th
component of $\Z^{m+1}$, we have found elements $r^0,\dots,r^m\in
\tilde{\sigma}^{\veee} \cap (\G\times\Z^{m+1})^{\ast}$.
\item[(iv)]
$\sigma\subseteq N_{\R}$ is the cone over $Q\subseteq \A$. Hence, the affine
embedding $\A \hookrightarrow \A\times\R^{m+1}\; (a\mapsto (a;\,1,\dots,1))$
induces an embedding of lattices $N \hookrightarrow \G\times\Z^{m+1}$ such
that $N= (\G\times\Z^{m+1}) \cap  \bigcap_{i=1}^m (r^i-r^0)^{\bot}$ and
$\sigma = \tilde{\sigma}\cap N_{\R}$.
\item[(v)]
$\{x^{r^1}-x^{r^0},\dots,x^{r^m}-x^{r^0}\}$ is a toric regular sequence in
$X:= \Spec \C[\tilde{\sigma}^{\veee}\cap (\G\times\Z^{m+1})^{\ast}]$,
and $Y$ is equal to the special fiber of the corresponding
(flat) map $X\rightarrow
\C^m$.
\end{itemize}
(The proof can be found in \S 4 of \cite{Mink}.)\\
\par

{\bf Remark:} (1)
The assumption that the degree $-\bar{r}$ has to be a primitive element of the
lattice $M$ is not essential. However, the description of the corresponding
toric regular sequences becomes slightly more complicated in the genaral case
(cf. \S 3 of \cite{Mink}) - and we do not need it in the present paper.\\
(2) The previous method yields deformations of degrees contained in
$-(\sigma^{\veee}\cap M)$ only. Nevertheless, $T^1$ can be non-trivial in
other degrees, too.\\
\par


\neu{32}
As a direct consequence we obtain\\
\par
{\bf Theorem:}
Let $Y$ be an affine toric Gorenstein variety induced
from a lattice polytope $Q$. Then, toric $m$-parameter
deformations of degree $-R^{\ast}$ correspond to Minkowski
decompositions of
$Q$ into a sum $Q=R_0+\dots +R_m$ of $m+1$ lattice polytopes.\\
The Kodaira-Spencer-map maps the parameter space $\C^m$ onto the linear
subspace
$\linh (\varrho(R_0),\dots ,\varrho(R_m)) \subseteq \tilde{T}^1(Q) =
T^1_Y(-R^{\ast})\subseteq T^1_Y$.\\
\par


{\bf Proof:}
$(\A,\G)$ defined in (i) of the previous theorem is exactly
that affine space containing our polytope $Q$. Moreover, $Q$ coincides with
the polyhedron $Q:= \sigma \cap \A$ defined in (ii).
Since $Q$ is a lattice polytope, the conditions for the summands $R_i$
(cf. (ii) of the previous theorem) are equivalent to the property of being
lattice polytopes, too.\\
Finally, the claim concerning the Kodaira-Spencer-map follows from the
definitions of the maps $\varrho$ and $\theta$ in \zitat{2}{5} and
\zitat{2}{6}, respectively, and from Theorem (5.2) of \cite{Mink}.
\hfill$\Box$\\
\par


\neu{33}
Let us focus on the special case of $\mbox{dim}\, Y=3$. Let $Y$ be given by
a 2-dimensional lattice polygon $Q=\mbox{Conv}(a^1,\dots,a^N)$ with
primitive edges
$\vec{v}_i:=a^{i+1}-a^i\;
(i\in \Z\left/_{\displaystyle \!\!N\,\Z}\right)$,
i.e. $Y$ has an isolated singularity in $0 \in Y$. Then, we obtain\\
\par
{\bf Theorem:} Non-trivial toric $m$-parameter deformations of $Y$
correspond to non-trivial Minkowski decompositions of $Q$ into a sum of
$m+1$ lattice polygons, i.e. to decompositions of the set of edges of $Q$
into a disjoint union of $m+1$ subsets each suming up to 0.\\
\par


{\bf Proof:} This is an immediate consequence of Theorem \zitat{3}{2} and
the fact that $T^1_Y$ is concentrated in degree $-R^\ast$ (cf.
\zitat{2}{9}). Nevertheless, beeing a little more carefully, for this
conclusion the following fact has to be used: Toric, regular sequences
inducing a trivial Kodaira-Spencer-map always yield trivial (standard)
deformations. This is proved in \S 6 of \cite{Mink}.
\hfill$\Box$\\
\par

%
%
\sect{Polarity and Examples}\label{s4}


\neu{41}
We start with recalling some general facts concerning the relation between
lattice polytopes and projective toric varieties (cf. Chapter 2 of
\cite{Oda}).\\
\par
Let $Q\subseteq (\A,\G)$ be a lattice polytope. Then, the inner normal fan
$\Sigma$ induces a projective toric variety ${\PP}(Q)$, and $Q$ itself
corresponds to an ample line bundle on it.\\
\par
Equivariant Weil divisors on
$\PP(Q)$ are described by maps $\Sigma^{(1)} \stackrel{h}{\rightarrow}\Z$.
Modulo principal divisors, they generate the whole divisor class group
$\mbox{Div}(\PP(Q))$. Let $D_h$ be a Weil divisor on $\PP(Q)$.
\begin{itemize}
\item[(i)]
$D_h$ is Cartier if and only if, on each top dimensional cone
$\alpha = \langle c^1,\dots,c^k\rangle \in\Sigma$, the map $h$ can be
represented as
\[
h(c^j) = \langle a_{\alpha},c^j\rangle \;(j=1,\dots,k) \quad \mbox{with }
a_{\alpha}\in\G.
\]
\item[(ii)]
A Cartier divisor $D_h$ is nef if and only if, moreover, $h\leq \langle
a_{\alpha},\bullet\rangle$ holds (for each top dimensional $\alpha\in\Sigma$)
on the whole 1-skeleton $\Sigma^{(1)}$. Then, the elements $a_{\alpha}$
form the vertex set of a polytope with inner normal cones containing the
corresponding cones of $\Sigma$.
\end{itemize}
On the other hand, if $Q'$ is a lattice polytope such that $\Sigma$ is a
subdivision of its inner normal fan (i.e. $Q'$ is a Minkowski summand of a
scalar multiple of $Q$), then we can use its vertices to define a map
$h(Q'):\Sigma^{(1)}\rightarrow\Z$ via (i). We obtain a nef Cartier divisor
$D_{Q'}$ on $\PP(Q)$ again.\\
The divisor $D_{Q'}$ is even ample if and only if $Q'$ and $Q$ induce the
same inner normal fan $\Sigma$. Equivalently, the elements $a_{\alpha}$ yield
different lattice points for different cones $\alpha$.\\
\par
Finally, we remark that variations by principal divisors correspond to
translations of the polytopes by lattice vectors only.\\
\par


\neu{42}
Let $Y$ be a 3-dimensional affine toric Gorenstein variety induced
by a lattice polygon $Q$, let
$Y$ having an isolated singularity in $0\in Y$. Then, we call
$\PP(Q)$ the polar variety of $Y$, it will be denoted by $\PP (Y)$.
\begin{itemize}
\item[(I)]
$\mbox{Pic}\,\PP (Y)$ equals the group generated by the lattice Minkowski
summands of scalar multiples of $Q$.
If we denote by
${\cal O}_{\PP (Y)}(1)$
the ample line bundle corresponding to $Q$ itself, then
\zitat{2}{9} tells us that
\[
T^1_Y \,=\; \mbox{Pic}\,{\PP}(Y) \left/_{\displaystyle \!\!{\cal O}(1)}
\right. \otimes_{\Z} \C.
\]
Let $P(Y)$ be the cone over $(\PP (Y),\,{\cal O}(1))$. The pull back of
${\cal O}(1)$  is a principal divisor on $P(Y)\setminus \{0\}$. Hence,
\[
T^1_Y \,=\; \mbox{Pic}(P(Y)\setminus \{0\}) \otimes_{\Z} \C.
\]
\item[(II)]
Theorem \zitat{3}{3} deals with Minkowski decompositions of $Q$ into a sum
of lattice polygons. In the language of the polar variety we obtain that
non-trivial toric $m$-parameter deformations of $Y$ correspond to non-trivial
decompositions of ${\cal O}_{\PP (Y)}(1)$ into a tensor product
\[
{\cal O}_{{\PP}(Y)}(1) = {\cal L}_0 \otimes \dots \otimes {\cal L}_m
\]
of $m+1$ nef invertible sheaves on $\PP (Y)$. The tangent plane to this
deformation inside the semi-universal base space $S$ is spanned by the
classes $[{\cal L}_0], \dots, [{\cal L}_m] \in \mbox{Pic}\,{\PP}(Y)
\left/_{\displaystyle \!\!{\cal O}(1)}\right. = \mbox{Pic}(P(Y)\setminus\{0\})
\subseteq T^1_Y$.
\vspace{2ex}
\end{itemize}
\par


{\bf Remark:}
The condition ``$0\in Y$ is an isolated singularity'' can be translated
into the
$\PP (Y)$-language, too: Each closed equivariant subvariety of
$\PP (Y)$
equals a linearly embedded projective space.\\
\par


\neu{43}
{\bf Conjecture:}
Take an arbitrary Minkowski decomposition of $Q$ into a sum of lattice
polytopes (equivalently: a decomposition of ${\cal O}_{\PP (Y)}(1)$ into a
tensor
product of nef line bundles), project the summands into $T^1_Y = \tilde{T}^1
(Q)\otimes_{\R}\C = \mbox{Pic}(P(Y)\setminus\{0\})\otimes_{\Z}\C$, and form
their linear
hull. Then, the union of all linear subspaces obtained in this way equals
the reduced base space of the semi-universal deformation of $Y$.\\
\par
The {\em cone} of Minkowski summands (i.e. the nef sheaves in $\mbox{Pic}\,
\PP (Y)$) contains much more information than the so-called {\em space}
of Minkowski summands of $Q$ (i.e.
$\mbox{Pic}(P(Y)\setminus\{0\})\otimes_{\Z}\R$).
Apart from the toric context, is there any such cone in deformation theory?
Projecting a certain interior point onto $0$ has to yield $T^1_Y$, then.\\
\par
Does the cone of Minkowski summands contain any information about the
non-reduced structure of the base space?\\
\par


\neu{44}
Finally, we want to present a special class of examples. We are looking for
those three-dimensional $Y$ that are, additionally to the usual assumptions,
cones over
projective toric varieties. (Do not mistake this property for $P(Y)$ being the
cone over $\PP (Y)$.)\\
In \S 4 of \cite{Bat} it is shown that these $Y$ can be characterized as the
cones over two-dimensional toric Fano varieties with Gorenstein singularities.
(Then, $P(Y)$ admits the same property.)\\
\par
The corresponding $Q$ are exactly
those lattice polygons containing one and only one interior lattice point
(``reflexive polygons''). Choosing this point as the origin, the polar
polygon $Q^{\veee}:= \{r\in \A^{\ast}|\,\langle Q,\,r\rangle \leq 1\}$ is a
lattice polygon, too. Then, $\stackrel{\veee}{\sigma}=\mbox{Cone}(Q^{\veee})$,
and
$Y$ is the cone over the projective variety corresponding to $Q^{\veee}$.\\
\par
Reflexive polygons were classified in (4.2) of \cite{Rob}. Our additional
assumption of
$Y$ having only an isolated singularity
causes that only five polygons
$Q$ survive from the original list (containing 16 items). Including the
polar polygons $Q^{\veee}$, we will see nine ones, however.\\
\par


{\bf (\ref{s4}.\ref{t44}.1)}\\
%
%
%
\unitlength=0.5mm
\linethickness{0.4pt}
\begin{picture}(51.00,51.00)(-60,-10)
\put(10.00,10.00){\circle*{1.00}}
\put(20.00,10.00){\circle*{1.00}}
\put(30.00,10.00){\circle*{1.00}}
\put(40.00,10.00){\circle*{1.00}}
\put(50.00,10.00){\circle*{1.00}}
\put(10.00,20.00){\circle*{1.00}}
\put(20.00,20.00){\circle*{1.00}}
\put(30.00,20.00){\circle*{1.00}}
\put(40.00,20.00){\circle*{1.00}}
\put(50.00,20.00){\circle*{1.00}}
\put(10.00,30.00){\circle*{1.00}}
\put(20.00,30.00){\circle*{1.00}}
\put(30.00,30.00){\circle*{1.00}}
\put(40.00,30.00){\circle*{1.00}}
\put(50.00,30.00){\circle*{1.00}}
\put(10.00,40.00){\circle*{1.00}}
\put(20.00,40.00){\circle*{1.00}}
\put(30.00,40.00){\circle*{1.00}}
\put(40.00,40.00){\circle*{1.00}}
\put(50.00,40.00){\circle*{1.00}}
\put(10.00,50.00){\circle*{1.00}}
\put(20.00,50.00){\circle*{1.00}}
\put(30.00,50.00){\circle*{1.00}}
\put(40.00,50.00){\circle*{1.00}}
\put(50.00,50.00){\circle*{1.00}}
\put(20.00,40.00){\line(0,-1){10.00}}
\put(20.00,30.00){\line(2,-1){20.00}}
\put(40.00,20.00){\line(0,1){10.00}}
\put(40.00,30.00){\line(-2,1){20.00}}
\put(30.00,0.00){\makebox(0,0)[cc]{Polygon $Q_1$}}
\end{picture}
%
%
%
\unitlength=0.5mm
\linethickness{0.4pt}
\begin{picture}(51.00,51.00)(-100,-10)
\put(10.00,10.00){\circle*{1.00}}
\put(20.00,10.00){\circle*{1.00}}
\put(30.00,10.00){\circle*{1.00}}
\put(40.00,10.00){\circle*{1.00}}
\put(50.00,10.00){\circle*{1.00}}
\put(10.00,20.00){\circle*{1.00}}
\put(20.00,20.00){\circle*{1.00}}
\put(30.00,20.00){\circle*{1.00}}
\put(40.00,20.00){\circle*{1.00}}
\put(50.00,20.00){\circle*{1.00}}
\put(10.00,30.00){\circle*{1.00}}
\put(20.00,30.00){\circle*{1.00}}
\put(30.00,30.00){\circle*{1.00}}
\put(40.00,30.00){\circle*{1.00}}
\put(50.00,30.00){\circle*{1.00}}
\put(10.00,40.00){\circle*{1.00}}
\put(20.00,40.00){\circle*{1.00}}
\put(30.00,40.00){\circle*{1.00}}
\put(40.00,40.00){\circle*{1.00}}
\put(50.00,40.00){\circle*{1.00}}
\put(10.00,50.00){\circle*{1.00}}
\put(20.00,50.00){\circle*{1.00}}
\put(30.00,50.00){\circle*{1.00}}
\put(40.00,50.00){\circle*{1.00}}
\put(50.00,50.00){\circle*{1.00}}
\put(20.00,40.00){\line(1,0){20.00}}
\put(40.00,40.00){\line(0,-1){20.00}}
\put(40.00,20.00){\line(-1,0){20.00}}
\put(20.00,20.00){\line(0,1){20.00}}
\put(30.00,0.00){\makebox(0,0)[cc]{Polar polygon $Q_1^{\veee}$}}
\end{picture}
\\
$Y_1$ is the cone over $\PP^1\times\PP^1$ embedded by ${\cal O}(2,2)$.\\
\par
$Q_1$ is a quadrangle, hence $T^1$ is one-dimensional. Moreover, $Q_1$ is the
Minkowski sum of two line segments, i.e. there is a really existing toric
1-parameter deformation. The total space is an isolated 4-dimensional cyclic
quotient singularity.\\
\par
In particular, the base space $S_1$ of the semi-universal
deformation of $Y_1$ equals $\C^1$.\\
\par


{\bf (\ref{s4}.\ref{t44}.2)}\\
%
%
%
\unitlength=0.5mm
\linethickness{0.4pt}
\begin{picture}(51.00,51.00)(-60,-10)
\put(10.00,10.00){\circle*{1.00}}
\put(20.00,10.00){\circle*{1.00}}
\put(30.00,10.00){\circle*{1.00}}
\put(40.00,10.00){\circle*{1.00}}
\put(50.00,10.00){\circle*{1.00}}
\put(10.00,20.00){\circle*{1.00}}
\put(20.00,20.00){\circle*{1.00}}
\put(30.00,20.00){\circle*{1.00}}
\put(40.00,20.00){\circle*{1.00}}
\put(50.00,20.00){\circle*{1.00}}
\put(10.00,30.00){\circle*{1.00}}
\put(20.00,30.00){\circle*{1.00}}
\put(30.00,30.00){\circle*{1.00}}
\put(40.00,30.00){\circle*{1.00}}
\put(50.00,30.00){\circle*{1.00}}
\put(10.00,40.00){\circle*{1.00}}
\put(20.00,40.00){\circle*{1.00}}
\put(30.00,40.00){\circle*{1.00}}
\put(40.00,40.00){\circle*{1.00}}
\put(50.00,40.00){\circle*{1.00}}
\put(10.00,50.00){\circle*{1.00}}
\put(20.00,50.00){\circle*{1.00}}
\put(30.00,50.00){\circle*{1.00}}
\put(40.00,50.00){\circle*{1.00}}
\put(50.00,50.00){\circle*{1.00}}
\put(20.00,20.00){\line(1,2){10.00}}
\put(30.00,40.00){\line(1,-1){10.00}}
\put(40.00,30.00){\line(-2,-1){20.00}}
\put(30.00,0.00){\makebox(0,0)[cc]{Polygon $Q_2$}}
\end{picture}
%
%
%
\unitlength=0.5mm
\linethickness{0.4pt}
\begin{picture}(51.00,51.00)(-100,-10)
\put(10.00,10.00){\circle*{1.00}}
\put(20.00,10.00){\circle*{1.00}}
\put(30.00,10.00){\circle*{1.00}}
\put(40.00,10.00){\circle*{1.00}}
\put(50.00,10.00){\circle*{1.00}}
\put(10.00,20.00){\circle*{1.00}}
\put(20.00,20.00){\circle*{1.00}}
\put(30.00,20.00){\circle*{1.00}}
\put(40.00,20.00){\circle*{1.00}}
\put(50.00,20.00){\circle*{1.00}}
\put(10.00,30.00){\circle*{1.00}}
\put(20.00,30.00){\circle*{1.00}}
\put(30.00,30.00){\circle*{1.00}}
\put(40.00,30.00){\circle*{1.00}}
\put(50.00,30.00){\circle*{1.00}}
\put(10.00,40.00){\circle*{1.00}}
\put(20.00,40.00){\circle*{1.00}}
\put(30.00,40.00){\circle*{1.00}}
\put(40.00,40.00){\circle*{1.00}}
\put(50.00,40.00){\circle*{1.00}}
\put(10.00,50.00){\circle*{1.00}}
\put(20.00,50.00){\circle*{1.00}}
\put(30.00,50.00){\circle*{1.00}}
\put(40.00,50.00){\circle*{1.00}}
\put(50.00,50.00){\circle*{1.00}}
\put(20.00,50.00){\line(0,-1){30.00}}
\put(20.00,20.00){\line(1,0){30.00}}
\put(50.00,20.00){\line(-1,1){30.00}}
\put(30.00,0.00){\makebox(0,0)[cc]{Polar polygon $Q_2^{\veee}$}}
\end{picture}
\\
$Y_2$ is the cone over $\PP^2$ embedded by ${\cal O}(3)$. Since $Q_2$ is a
triangle, $Y_2$ is rigid.\\
\par


{\bf (\ref{s4}.\ref{t44}.3)}\\
%
%
%
\unitlength=0.5mm
\linethickness{0.4pt}
\begin{picture}(51.00,51.00)(-60,-10)
\put(10.00,10.00){\circle*{1.00}}
\put(20.00,10.00){\circle*{1.00}}
\put(30.00,10.00){\circle*{1.00}}
\put(40.00,10.00){\circle*{1.00}}
\put(50.00,10.00){\circle*{1.00}}
\put(10.00,20.00){\circle*{1.00}}
\put(20.00,20.00){\circle*{1.00}}
\put(30.00,20.00){\circle*{1.00}}
\put(40.00,20.00){\circle*{1.00}}
\put(50.00,20.00){\circle*{1.00}}
\put(10.00,30.00){\circle*{1.00}}
\put(20.00,30.00){\circle*{1.00}}
\put(30.00,30.00){\circle*{1.00}}
\put(40.00,30.00){\circle*{1.00}}
\put(50.00,30.00){\circle*{1.00}}
\put(10.00,40.00){\circle*{1.00}}
\put(20.00,40.00){\circle*{1.00}}
\put(30.00,40.00){\circle*{1.00}}
\put(40.00,40.00){\circle*{1.00}}
\put(50.00,40.00){\circle*{1.00}}
\put(10.00,50.00){\circle*{1.00}}
\put(20.00,50.00){\circle*{1.00}}
\put(30.00,50.00){\circle*{1.00}}
\put(40.00,50.00){\circle*{1.00}}
\put(50.00,50.00){\circle*{1.00}}
\put(20.00,20.00){\line(0,1){10.00}}
\put(20.00,30.00){\line(2,1){20.00}}
\put(40.00,40.00){\line(-1,-2){10.00}}
\put(30.00,20.00){\line(-1,0){10.00}}
\put(30.00,0.00){\makebox(0,0)[cc]{Polygon $Q_3$}}
\end{picture}
%
%
%
\unitlength=0.5mm
\linethickness{0.4pt}
\begin{picture}(51.00,51.00)(-100,-10)
\put(10.00,10.00){\circle*{1.00}}
\put(20.00,10.00){\circle*{1.00}}
\put(30.00,10.00){\circle*{1.00}}
\put(40.00,10.00){\circle*{1.00}}
\put(50.00,10.00){\circle*{1.00}}
\put(10.00,20.00){\circle*{1.00}}
\put(20.00,20.00){\circle*{1.00}}
\put(30.00,20.00){\circle*{1.00}}
\put(40.00,20.00){\circle*{1.00}}
\put(50.00,20.00){\circle*{1.00}}
\put(10.00,30.00){\circle*{1.00}}
\put(20.00,30.00){\circle*{1.00}}
\put(30.00,30.00){\circle*{1.00}}
\put(40.00,30.00){\circle*{1.00}}
\put(50.00,30.00){\circle*{1.00}}
\put(10.00,40.00){\circle*{1.00}}
\put(20.00,40.00){\circle*{1.00}}
\put(30.00,40.00){\circle*{1.00}}
\put(40.00,40.00){\circle*{1.00}}
\put(50.00,40.00){\circle*{1.00}}
\put(10.00,50.00){\circle*{1.00}}
\put(20.00,50.00){\circle*{1.00}}
\put(30.00,50.00){\circle*{1.00}}
\put(40.00,50.00){\circle*{1.00}}
\put(50.00,50.00){\circle*{1.00}}
\put(20.00,20.00){\line(1,0){20.00}}
\put(40.00,20.00){\line(0,1){30.00}}
\put(40.00,50.00){\line(-1,-1){20.00}}
\put(20.00,30.00){\line(0,-1){10.00}}
\put(30.00,0.00){\makebox(0,0)[cc]{Polar polygon $Q_3^{\veee}$}}
\end{picture}
\\
$Y_3$ is the cone over the Del Pezzo surface of degree 8 (the blowing up of
$(\PP^2,{\cal O}(3))$ in one point).\\
\par
The vector space $T^1$ is one-dimensional. The two-dimensional cone of the
rational Minkowski summands of scalar multiples of $Q_3$ is generated by two
triangles.\\
\par
However, there are no lattice polygons that are non-trivial Minkowski summands
of $Q_3$.
That means, $Y_3$ does not admit any toric deformation at all.\\
Indeed, as Duco van Straten has computed with Macaulay, the semi-universal
base space $S_3$ of $Y_3$ equals
$\Spec \C[\varepsilon]/_{\displaystyle \varepsilon^2}$.\\
\par


{\bf (\ref{s4}.\ref{t44}.4)}\\
%
%
%
\unitlength=0.5mm
\linethickness{0.4pt}
\begin{picture}(51.00,51.00)(-60,-10)
\put(10.00,10.00){\circle*{1.00}}
\put(20.00,10.00){\circle*{1.00}}
\put(30.00,10.00){\circle*{1.00}}
\put(40.00,10.00){\circle*{1.00}}
\put(50.00,10.00){\circle*{1.00}}
\put(10.00,20.00){\circle*{1.00}}
\put(20.00,20.00){\circle*{1.00}}
\put(30.00,20.00){\circle*{1.00}}
\put(40.00,20.00){\circle*{1.00}}
\put(50.00,20.00){\circle*{1.00}}
\put(10.00,30.00){\circle*{1.00}}
\put(20.00,30.00){\circle*{1.00}}
\put(30.00,30.00){\circle*{1.00}}
\put(40.00,30.00){\circle*{1.00}}
\put(50.00,30.00){\circle*{1.00}}
\put(10.00,40.00){\circle*{1.00}}
\put(20.00,40.00){\circle*{1.00}}
\put(30.00,40.00){\circle*{1.00}}
\put(40.00,40.00){\circle*{1.00}}
\put(50.00,40.00){\circle*{1.00}}
\put(10.00,50.00){\circle*{1.00}}
\put(20.00,50.00){\circle*{1.00}}
\put(30.00,50.00){\circle*{1.00}}
\put(40.00,50.00){\circle*{1.00}}
\put(50.00,50.00){\circle*{1.00}}
\put(20.00,30.00){\line(0,1){10.00}}
\put(20.00,40.00){\line(1,0){10.00}}
\put(30.00,40.00){\line(1,-1){10.00}}
\put(40.00,30.00){\line(-1,-1){10.00}}
\put(30.00,20.00){\line(-1,1){10.00}}
\put(30.00,0.00){\makebox(0,0)[cc]{Polygon $Q_4$}}
\end{picture}
%
%
%
\unitlength=0.5mm
\linethickness{0.4pt}
\begin{picture}(51.00,51.00)(-100,-10)
\put(10.00,10.00){\circle*{1.00}}
\put(20.00,10.00){\circle*{1.00}}
\put(30.00,10.00){\circle*{1.00}}
\put(40.00,10.00){\circle*{1.00}}
\put(50.00,10.00){\circle*{1.00}}
\put(10.00,20.00){\circle*{1.00}}
\put(20.00,20.00){\circle*{1.00}}
\put(30.00,20.00){\circle*{1.00}}
\put(40.00,20.00){\circle*{1.00}}
\put(50.00,20.00){\circle*{1.00}}
\put(10.00,30.00){\circle*{1.00}}
\put(20.00,30.00){\circle*{1.00}}
\put(30.00,30.00){\circle*{1.00}}
\put(40.00,30.00){\circle*{1.00}}
\put(50.00,30.00){\circle*{1.00}}
\put(10.00,40.00){\circle*{1.00}}
\put(20.00,40.00){\circle*{1.00}}
\put(30.00,40.00){\circle*{1.00}}
\put(40.00,40.00){\circle*{1.00}}
\put(50.00,40.00){\circle*{1.00}}
\put(10.00,50.00){\circle*{1.00}}
\put(20.00,50.00){\circle*{1.00}}
\put(30.00,50.00){\circle*{1.00}}
\put(40.00,50.00){\circle*{1.00}}
\put(50.00,50.00){\circle*{1.00}}
\put(20.00,30.00){\line(0,-1){10.00}}
\put(20.00,20.00){\line(1,0){20.00}}
\put(40.00,20.00){\line(0,1){20.00}}
\put(40.00,40.00){\line(-1,0){10.00}}
\put(30.00,40.00){\line(-1,-1){10.00}}
\put(30.00,0.00){\makebox(0,0)[cc]{Polar polygon $Q_4^{\veee}$}}
\end{picture}
\\
$Y_4$ is the cone over the Del Pezzo surface of degree 7 (obtained from
$(\PP^2,{\cal O}(3))$ by blowing up two points, or from $(\PP^1\times\PP^1,
{\cal O}(2,2))$ by blowing up one point).\\
\par
$T^1$ is two-dimensional, but $Q_4$ admits one decomposition into a
Minkowski sum of two lattice polygons only. $Q_4$ equals the sum of a line
segment and a triangle - this yields a 1-parameter deformation of $Y_4$,
the total space is the cone over $\PP ({\cal O}_{\PP^2}\oplus {\cal O}
_{\PP^2}(1))$.\\
\par
The semi-universal base space $S_4$ is a complex line with one embedded
component (computed by Duco van Straten using Macaulay).\\
\par


{\bf (\ref{s4}.\ref{t44}.5)}\\
%
%
%
\unitlength=0.5mm
\linethickness{0.4pt}
\begin{picture}(51.00,51.00)(-90,-10)
\put(10.00,10.00){\circle*{1.00}}
\put(20.00,10.00){\circle*{1.00}}
\put(30.00,10.00){\circle*{1.00}}
\put(40.00,10.00){\circle*{1.00}}
\put(50.00,10.00){\circle*{1.00}}
\put(10.00,20.00){\circle*{1.00}}
\put(20.00,20.00){\circle*{1.00}}
\put(30.00,20.00){\circle*{1.00}}
\put(40.00,20.00){\circle*{1.00}}
\put(50.00,20.00){\circle*{1.00}}
\put(10.00,30.00){\circle*{1.00}}
\put(20.00,30.00){\circle*{1.00}}
\put(30.00,30.00){\circle*{1.00}}
\put(40.00,30.00){\circle*{1.00}}
\put(50.00,30.00){\circle*{1.00}}
\put(10.00,40.00){\circle*{1.00}}
\put(20.00,40.00){\circle*{1.00}}
\put(30.00,40.00){\circle*{1.00}}
\put(40.00,40.00){\circle*{1.00}}
\put(50.00,40.00){\circle*{1.00}}
\put(10.00,50.00){\circle*{1.00}}
\put(20.00,50.00){\circle*{1.00}}
\put(30.00,50.00){\circle*{1.00}}
\put(40.00,50.00){\circle*{1.00}}
\put(50.00,50.00){\circle*{1.00}}
\put(20.00,20.00){\line(1,0){10.00}}
\put(30.00,20.00){\line(1,1){10.00}}
\put(40.00,30.00){\line(0,1){10.00}}
\put(40.00,40.00){\line(-1,0){10.00}}
\put(30.00,40.00){\line(-1,-1){10.00}}
\put(20.00,30.00){\line(0,-1){10.00}}
\put(30.00,0.00){\makebox(0,0)[cc]{Polygon $Q_5=Q_5^{\veee}$}}
\end{picture}
\\
$Y_5$ is the cone over the Del Pezzo surface of degree 6 (obtained by blowing
up the projective variety of (\ref{s4}.\ref{t44}.4) in one more point).\\
$T^1$ is three-dimensional, and $Q_5$ admits two different extremal Minkowski
decompositions:
\begin{itemize}
\item[(i)]
$Q_5$ equals the sum of two triangles, the corresponding 1-parameter family
admits the cone over $\PP^1\times\PP^1\times\PP^1$ as its total space.
\item[(ii)]
$Q_5$ also equals the sum of three line segments. This corresponds to a
two-parameter family with the cone over $\PP^2\times\PP^2$ as its total space.
\end{itemize}
Again, Duco van Straten has computed the semi-universal base space - it is
reduced and equals the
transversal union of a complex plane with a complex line. These components
correspond to the toric deformations we have already seen.

%
%


\begin{thebibliography}{2cm}

\bibitem[Al 1]{T1} Altmann, K.: Computation of the vector space $T^1$ for
affine toric varieties.\\
J. Pure Appl. Algebra (to appear).

\bibitem[Al 2]{Mink} Altmann, K.: Minkowski sums and homogeneous deformations
of toric varieties.\\
Preprint 93-1, Humboldt-Universit\"at Berlin or\\
Preprint N\underline{o}. 22, Europ\"aisches Singularit\"atenprojekt,
Berlin 1993.

\bibitem[Ar]{Arndt} Arndt, J.: Verselle Deformationen zyklischer
Quotientensingularit\"aten. \\
Dissertation, Universit\"at Hamburg, 1988.

\bibitem[Ba]{Bat} Batyrev, V.: Dual polyhedra and the mirror symmetry for
Calabi-Yau hypersurfaces in toric varieties.\\
Preprint, University of Essen, 1992.

\bibitem[Ch]{Ch} Christophersen, J.A.: Obstruction spaces for rational
singularities and deformations of cyclic quotients. \\
Thesis, University of Oslo, 1989/90.

\bibitem[Ke]{Ke} Kempf, G., Knudsen, F., Mumford, D., Saint-Donat, B.:
Toroidal Embeddings I. \\
Lecture Notes in Mathematics {\bf 339}, Springer-Verlag,
Berlin-Heidelberg-New York, 1973.

\bibitem[Ko]{Rob} Koelman, R.J.: The number of moduli of families of curves on
toric surfaces.\\
Proefschrift, Nijmegen, 1991.

\bibitem[Od]{Oda} Oda, T.: Convex bodies and algebraic geometry.\\
Ergebnisse der Mathematik und ihrer Grenzgebiete (3/15), Springer-Verlag, 1988.

\bibitem[Ri]{Riem} Riemenschneider, O.: Deformationen von
Quotienten\-singularit\"aten (nach zyklischen Gruppen).\\
Math. Ann. {\bf 209} (1974), 211-248.

\bibitem[Sm]{Sm} Smilansky, Z.: Decomposability of Polytopes and Polyhedra.\\
Geometriae Dedicata. {\bf 24} (1987), 29-49.

\bibitem[St]{St} Stevens, J.: On the versal deformation of cyclic quotient
singularities. \\
Preprint Hamburg.



\end{thebibliography}
\end{document}